\documentclass[10pt,english, a4paper]{article}
\usepackage {amssymb}
\usepackage[intlimits] {amsmath}
\usepackage {babel}
\usepackage {exscale}
\usepackage{epsfig}
\usepackage{graphicx}
\usepackage {theorem}
\usepackage[latin1]{inputenc}
\usepackage{makeidx}
\theoremstyle{break}

\renewcommand{\v}{\mathbf{v}}

\renewcommand{\epsilon}{\varepsilon}

\newcommand{\norm}[1]{\|#1\|}

\newtheorem{lem}{Lemma}[section]

\newtheorem{rem}[lem]{Remark}
\newtheorem{defi}[lem]{Definition}
\newtheorem{prop}[lem]{Proposition}

\newtheorem{theo}[lem]{Theorem}
\newtheorem{conj}[lem]{Conjecture}

\newtheorem{cor}[lem]{Corollary}
\newtheorem{exa}[lem]{Example}

\newcommand{\leref}[1]{Lemma \ref{#1}}
\newcommand{\theref}[1]{Theorem \ref{#1}}
\newcommand{\conjref}[1]{Conjecture \ref{#1}}

\newcommand{\coref}[1]{Corollary \ref{#1}}

\newcommand{\propref}[1]{Proposition \ref{#1}}

\newcommand{\deref}[1]{Definition \ref{#1}}

\newcommand{\R}{\mathbb{R}}

\newcommand{\N}{\mathbb{N}}

\newcommand{\C}{\mathbb{C}}

\DeclareMathOperator{\Log}{Log}

\DeclareMathOperator{\GL}{GL}
\DeclareMathOperator{\SO}{SO}

\DeclareMathOperator{\diag}{diag}
\DeclareMathOperator{\sym}{sym}

\DeclareMathOperator{\gl}{\mathfrak{gl}}

\DeclareMathOperator{\polar}{polar}
\newcommand{\qed}{\hfill $\blacksquare$\\}

\DeclareMathOperator{\Det}{det}

\renewcommand{\det}[1]{ {\Det[{#1}]} }

\DeclareMathOperator{\dist}{dist}

\theoremstyle{break}
\allowdisplaybreaks[1]
\makeatletter
\@addtoreset{equation}{section}

\makeindex

\begin{document}
\title{On the generalized sum of squared logarithms inequality}
 \author{Waldemar Pompe%
 \thanks{Corresponding author: Waldemar Pompe, Institute of Mathematics, University of Warsaw, ul. Banacha 2, 02-097 Warszawa, Poland, email: pompe@mimuw.edu.pl.} 
\, and 
Patrizio Neff%
\thanks {Patrizio Neff, Head of Chair for Nonlinear Analysis and Modelling, Fakult\"at f\"ur Mathematik, Universität Duisburg-Essen, Campus Essen, Thea-Leymann Str. 9, 45127 Essen, Germany, email: patrizio.neff@uni-due.de, Phone +49 201 183 4243, Fax: +49 201 183 4394}
 }
\maketitle

\begin{abstract}
Assume $n\geq 2$. Consider the elementary symmetric polynomials $e_k(y_1,y_2,\ldots, y_n)$ and denote by 
$E_0,E_1,\ldots,E_{n-1}$ the elementary symmetric polynomials in reverse order
\begin{align*}
    E_k(y_1,y_2,\ldots,y_n):=e_{n-k}(y_1,y_2,\ldots,y_n)=\hspace*{-3.5mm}\sum_{i_1<\ldots<i_{n-k}} y_{i_1}y_{i_2}\ldots y_{i_{n-k}}\, ,
\quad k\in \{0,1,\ldots,n{-}1     \}\, .
\end{align*}
Let moreover $S$ be a nonempty subset of $\{0,1,\ldots,n{-}1\}$.
We investigate necessary and sufficient conditions on the function $f\colon\,I\to\R$,
where $I\subset\R$ is an interval, such that the inequality 
\begin{align}
\label{abstract_inequality}
   f(a_1)+f(a_2)+\ldots+f(a_n)\leq f(b_1)+f(b_2)+\ldots+f(b_n) \tag{*}
\end{align}   
holds for all $a=(a_1,a_2,\ldots,a_n)\in I^n$ and
$b=(b_1,b_2,\ldots,b_n)\in I^n$ satisfying
$$E_k(a)< E_k(b) \ \hbox{for } k\in S\quad \hbox{and} \quad E_k(a)=E_k(b) 
\quad\hbox{for } k\in \{0,1,\ldots,n{-}1     \}\setminus S\, .$$
As a corollary, we obtain \eqref{abstract_inequality} if $2\leq n\leq 4$,
$f(x)=\log^2x$ and $S=\{1,\dotsc,n-1\}$, which is the sum of squared logarithms inequality previously known for $2\le n\le 3$.
\end{abstract}
\vspace{1cm}


\tableofcontents

\section{Introduction - the sum of squared logarithms inequality}
In a previous contribution \cite{Neff_log_inequality13} the sum of squared logarithms inequality has been introduced and proved for the particular cases $n=2,3$. For $n=3$ it reads: let $a_1,a_2,a_3,b_1,b_2,b_3>0$ be given positive numbers such that
\begin{align}
   a_1+a_2+a_3 &\le b_1+b_2+b_3  \, , \notag\\
   a_1\, a_2+a_1\, a_3+ a_2\, a_3 & \le b_1\, b_2+b_1\, b_3+ b_2\, b_3 \, ,\notag \\
    a_1 \, a_2\, a_3&=b_1\, b_2\, b_3 \, . \notag
\end{align}
Then
\begin{align}
 \log^2 a_1+\log^2 a_2+\log^2 a_3 \le  \log^2 b_1+\log^2 b_2+\log^2 b_3   \, .\notag
 \end{align}
The general form of this inequality can be conjectured as follows. 
\begin{defi}
	The \emph{standard elementary symmetric polynomials} $e_1,\ldots,e_{n-1}, e_n$ are
	\begin{align}
	e_k(y_1,\ldots,y_n)=\sum_{1\le j_1<j_2<\ldots<j_k\le n}      y_{j_1}\cdot y_{j_2}\ldots \cdot y_{j_k}\, ,
	\quad k\in\{1,2,\ldots,n\}\,;
	\end{align}
	note that $e_n=y_1\cdot y_2\ldots \cdot y_n$.
\end{defi}
\begin{conj}[Sum of squared logarithms inequality]
 \label{log_inequality}
 Let $a_1,a_2,\ldots,a_n$, $b_1,b_2,\ldots,b_n$ be given positive numbers. Then the condition
\[
  e_k(a_1,\ldots, a_n) \le e_k(b_1,\ldots, b_n), \quad k\in \{1,2,\ldots,n-1\}, \quad e_n(a_1,\ldots, a_n) = e_n(b_1,\ldots, b_n)  \notag\\
\]
implies that
\[
    \sum_{i=1}^n \log^2 a_i \le  \sum_{i=1}^n \log^2 b_i\, .
\]
\end{conj}

\begin{rem}
 \label{injectivity_remark}
Note that the conclusions of \conjref{log_inequality} are trivial provided we have equality everywhere, i.e.
\begin{align}
  e_k(a_1,\ldots, a_n) &= e_k(b_1,\ldots, b_n), \quad k\in \{1,2,\ldots,n\}\, .
  \end{align}
In this case, the coefficients $a_1,\ldots a_n, b_1,\ldots b_n$ are equal up to permutations, which can be seen by looking at the 
characteristic polynomials of two matrices with eigenvalues $a_1,\dotsc,a_n$ and $b_1,\dotsc,b_n$. From this perspective, having equality just in the last product $e_n$ and strict inequality else seems to be the most difficult case.
\end{rem}

Based on extensive random sampling on $\R_+^n$ for small numbers $n$ it has been conjectured that \conjref{log_inequality} might be true for arbitrary $n\in\N$. The sum of squared logarithms inequality has immediate important applications in matrix analysis (\cite{Neff_Nakatsukasa_logpolar13}, see also \cite{Lankeit_Neff_Nakatsukasa_logpolar13b}) as well as in nonlinear elasticity theory \cite{Neff_Ghiba_ZAMP14,Neff_Lankeit_Ghiba_JEL14,Neff_Grioli14,Neff_Osterbrink_hencky13}. In matrix analysis it implies that the global minimizer over all rotations to
\begin{align}
 \inf_{Q\in\SO(n)} \norm{\sym_*\Log Q^T \, F}^2=\norm{\sqrt{F^TF}}^2
\end{align}
at given $F\in\GL^+(n)$ is realized by the orthogonal factor $R=\polar(F)$ (such that $R^T\,F=\sqrt{F^TF}$). Here, $\norm{X}^2:=\sum_{i,j=1}^n X_{ij}^2$ denotes the Frobenius matrix norm and $\Log: \GL(n) \to \gl(n)=\R^{n\times n}$ is the multivalued matrix-logarithm, i.e. any solution $Z=\Log X\in \C^{n\times n}$ of $\exp(Z)=X$ and $\sym_*(Z)=\frac{1}{2}\left( Z^*+Z\right)$.

Recently, the case $n=2$ was used to verify the polyconvexity condition in nonlinear elasticity \cite{Neff_Lankeit_Ghiba_JEL14,Neff_Ghiba_ZAMP14} for a certain class of isotropic energy functions. For more background information on the sum of squared logarithms inequality we refer the reader to \cite{Neff_log_inequality13}.\\

In this paper we extend the investigation as to the validity of \conjref{log_inequality} by considering arbitrary functions $f$ instead of $f(x)=\log^2 x$. We formulate this more general problem and we are able to extend \conjref{log_inequality} to the case $n=4$. The same methods should also be useful for 
proving the statement for $n=5,6$. However, the necessary technicalities prevent us from discussing these cases in this paper.

In addition, we present ideas which might be helpful in attacking the fully general case, namely arbitrary $f$ and arbitrary $n$. 

\section{The generalized inequality}

\def\C{{\cal C}}
\def\inter{{\rm int}\,}
\def\ve{\varepsilon}
\def\diag{{\rm diag}\,}
\def\dist{\hbox{\rm dist}\,}
\def\norm#1{|\!|#1|\!|}
\def\Om{\Omega}
\def\v{\vbox to13pt{}}

In order to generalize \conjref{log_inequality} 
in the directions hinted at in the introduction, 
we consider from now on a non-standard definition of 
the elementary symmetric polynomials. In fact, for $n\geq 2$ 
it will be more convenient for us to reverse their numbering and define $E_0,E_1,\ldots,E_{n-1}$ by
\begin{align}
  \label{non_standard_polynomials}
     E_k(y_1,\ldots y_n):=e_{n-k}( y_1,\ldots, y_n)=\sum_{i_1<\ldots<i_{n-k}} y_{i_1}\cdot y_{i_2}\ldots \cdot y_{i_{n-k}}\, ,
\quad k\in\{0,1,\ldots,{n-1}\}\,.
\end{align}
In particular, now
\begin{align}
      E_0(y_1,\ldots, y_n)&:=e_{n}( y_1,\ldots, y_n)=y_1\cdot y_2\cdot\ldots \cdot y_n\, ,\notag\\
      E_{n-1}(y_1,\ldots, y_n)&:=e_{1}( y_1,\ldots, y_n)=y_1+y_2+\ldots+y_n\, .
\end{align}

\noindent
Let $I\subset \R$ be an open interval and let
\begin{align}
    \Delta_n:=\{y=(y_1,y_2,\ldots,y_n)\in I^n \,|\, y_1\leq y_2\leq\ldots\leq y_n\}\,.
\end{align}    
Let $S$ be a nonempty subset of $\{0,1,\ldots,n{-}1\}$ and assume that
$a,b\in\Delta_n$ are such that
\begin{align}
\label{inequality_1}
    E_k(a)< E_k(b) \quad\hbox{for } k\in S\qquad \hbox{and} \qquad E_k(a)=E_k(b) 
\quad\hbox{for } k\in \{0,1,\ldots,n{-}1\}\setminus S\,.
\end{align}
In this section we investigate necessary and sufficient conditions for a 
(smooth) function $f\colon\,I\to\R$, such that the inequality
$$f(a_1)+f(a_2)+\ldots+f(a_n)\leq f(b_1)+f(b_2)+\ldots+f(b_n)$$
holds for all $a,b\in\Delta_n$ satisfying assumption \eqref{inequality_1}.

\begin{rem}
The formulation of the above problem 
has a certain monotonicity structure: we assume that ``$E(a)<E(b)$''
and want to prove that ``$F(a)<F(b)$''. Therefore our idea is to consider
a curve $y$ connecting the points $a$ and $b$, such that $E(y(t))$ ``increases''.
Then the function $g(t)=F(y(t))$ should also increase and therefore 
$g'(t)>0$ must hold. From this we are able to derive necessary and 
sufficient conditions on the function $f$.
\end{rem}

This approach motivates the following definition.
\begin{defi}[$b$ \textnormal{\textit{dominates}} $a$, $a\preceq b$]
\label{definition_1}
Let $a,b\in\Delta_n$. We will say that $b$ {\it dominates\/} $a$, and denote $a\preceq b$, if there exists
a piecewise differentiable mapping $y\colon\,[0,1]\to\Delta_n$ 
(i.e. $y$ is continuous on $[0,1]$ and differentiable in all but at most countably many points)
such that
$y(0)=a$, $y(1)=b$, $y_i(t)\neq y_j(t)$ for $i\neq j$ and all but at most countably many $t\in[0,1]$ 
and the functions
$$A_k(t):=E_k(y(t))\, , \qquad k\in \{0,1,\ldots,n{-}1\}  $$
are non-decreasing on the interval $[0,1]$. 
\end{defi}

If $a\preceq b$, then $E_k(a)=A_k(0)\leq A_k(1)=E_k(b)$,
so it follows from \deref{definition_1} that 
$a, b$ satisfy assumption \eqref{inequality_1} with $S$ being the set of all $k$ for which
$A_k(t)$ is not a constant function on $[0,1]$.

\medskip
We are ready to formulate the main results of this section.

\begin{theo}
\label{theorem_1}
Assume that $a,b\in \Delta_n$ and let $a\preceq b$.  
Let $S\subseteq\{0,1,\ldots,n{-}1\}$ 
denote the set of all integers $k$ with $E_k(a)<E_k(b)$.
Moreover, assume that $f\in C^n(I)$ be such that 
\begin{align}
\label{differential_inequality_2}
     (-1)^{n+k}(x^k f'(x))^{(n-1)}\leq 0 \quad \hbox{for all }x\in I \hbox{ and all } k\in S\,.
\end{align}
Then the following inequality holds:
\begin{align}
\label{inequality_3} 
      f(a_1)+f(a_2)+\ldots+f(a_n)\leq f(b_1)+f(b_2)+\ldots+f(b_n)\,.
\end{align}
\end{theo}

A partially reverse statement is also true.
\begin{theo}
\label{theorem_2}
Let $f\in C^n(I)$ be such that the inequality 
\begin{align} 
      f(a_1)+f(a_2)+\ldots+f(a_n)\leq f(b_1)+f(b_2)+\ldots+f(b_n)
\end{align}
holds for all $a,b\in \Delta_n$ satisfying 
\begin{align}
\label{dupa}
    E_k(a)\leq E_k(b) \quad\hbox{for } k\in S\qquad \hbox{and} \qquad E_k(a)=E_k(b) 
\quad\hbox{for } k\in \{0,1,\ldots,n{-}1\}\setminus S
\end{align}
for some subset $S\subseteq\{0,1,\ldots,n{-}1\}$. 
Then $f$ satisfies property \eqref{differential_inequality_2}, i.e.
\begin{align}
     (-1)^{n+k}(x^k f'(x))^{(n-1)}\leq 0 \quad \hbox{for all }x\in I \hbox{ and all } k\in S\,.
\end{align}
\end{theo}

In this respect, we can formulate another conjecture:
\begin{conj}
\label{missing_conjecture}
Let $S$ be a nonempty subset of $\{0,1,\ldots,n{-}1\}$ and assume that
$a,b\in\Delta_n$ are such that \eqref{inequality_1} is satisfied, i.e.
\begin{align*}
    E_k(a)< E_k(b) \quad\hbox{for } k\in S\qquad \hbox{and} \qquad E_k(a)=E_k(b) 
\quad\hbox{for } k\in \{0,1,\ldots,n{-}1\}\setminus S\,.
\end{align*}
Then there exists a curve $y$ satisfying the conditions from \deref{definition_1} and thus $a\preceq b$.
\end{conj}

\begin{rem}
 In concrete applications of \theref{theorem_1} and \theref{theorem_2}
one would like to know whether condition \eqref{inequality_1} already implies $a\preceq b$. %
This is \conjref{missing_conjecture}. Unfortunately, we are able to prove \conjref{missing_conjecture} only for $2\leq n\leq 4$, $I=(0,\infty)$ and
$S\subseteq\{1,2,\ldots,n{-}1\}$ (see the next section). 
\end{rem}

\smallskip
\begin{exa}
It is easy to see that if $I=(0,\infty)$ then the function $f(x)=\log^2x$
satisfies property \eqref{differential_inequality_2} for $S=\{1,2,\ldots,n{-}1\}$. Indeed, we proceed by induction on $n$.
For $n=2$ and $k=1$ the property is immediate. Moreover, for $k\geq2$ and $n\geq3$ we get
\begin{align}
   (-1)^{n+k}(x^k f'(x))^{(n-1)}&=2(-1)^{n+k}(x^{k-1}\log x)^{(n-1)}  \\
   &=2(-1)^{n+k}((k-1)x^{k-2}\log x)^{(n-2)}+2(-1)^{n+k}(x^{k-2})^{(n-2)}\leq 0  \notag
   \end{align}
by the induction hypothesis, since the second summand vanishes. 
It remains to check property \eqref{differential_inequality_2} for $k=1$, which
is also immediate.

Note also that property \eqref{differential_inequality_2} is not true for $k=0$. Therefore \theref{theorem_1} and \theref{theorem_2} 
for $f(x)=\log^2x$ attain the following formulation:
\end{exa}

\begin{cor}
\label{corollary_1}
Assume that $a,b\in \R_+^n$ be such that $a\preceq b$ and
$a_1a_2\ldots a_n=b_1b_2\ldots b_n\,.$ Then
$$\log^2(a_1)+\log^2(a_2)+\ldots+\log^2(a_n)\leq 
\log^2(b_1)+\log^2(b_2)+\ldots+\log^2(b_n)\,$$
and this inequality fails, if the constraint 
$a_1a_2\ldots a_n=b_1b_2\ldots b_n$ is replaced by
the weaker one $a_1a_2\ldots a_n\leq b_1b_2\ldots b_n\,.$
\end{cor}

In order to see that the weaker condition is not sufficient for the inequality to hold, consider the case
\[
	a = (\tfrac1n,\dotsc,\tfrac1n)\,, \qquad b=(1,\dotsc,1)\,.
\]
Then $a \preceq b$ and $a_1a_2\ldots a_n\leq b_1b_2\ldots b_n$, but
\[
	\log^2(a_1)+\log^2(a_2)+\ldots+\log^2(a_n) = n\,\log^2(n) > 0 = \log^2(b_1)+\log^2(b_2)+\ldots+\log^2(b_n)\,.
\]

\begin{rem}
Corollary \ref{corollary_1} is a weaker statement than \conjref{log_inequality} since we assume that $a\preceq b$. If \conjref{missing_conjecture} is true, then \conjref{log_inequality} follows.
\end{rem}

\smallskip
\begin{exa}
The function $f(x)=x^p$ $(x>0)$ with $p\in(0,1)$
satisfies property \eqref{differential_inequality_2} for the set $S=\{0,1,\ldots,n{-}1\}$.
Indeed, for each $n\geq2$ and $0\leq k\leq n-1$, we have
$$(-1)^{n+k}(x^kf'(x))^{(n-1)}=
(-1)^{n+k}p(k+p-1)(k+p-2)\ldots(k+p-(n{-}1))x^{k+p-n}\,.$$
The above product is not greater than $0$,
because among the factors
$k+p-1,k+p-2,\ldots,k+p-(n{-}1)$ there are exactly
$n-1-k$ negative ones. 

Similarly, the function $f(x)=x^p$ for $p\in(-1,0)$
satisfies property \eqref{differential_inequality_2} for the set $S=\{1,2,\ldots,{n{-}1}\}$,
because $p<0$ and among the factors
$k+p-1,k+p-2,\ldots,k+p-(n{-}1)$ there are exactly
$n-k$ negative ones. On the other hand, property \eqref{differential_inequality_2} is not true for $k=0$.\\

Thus, similarly like above, we have
\end{exa}

\begin{cor}
\label{corollary_2}
Assume that $a,b\in (0,\infty)^n$ be such that $a\preceq b$ and
$a_1a_2\ldots a_n=b_1b_2\ldots b_n\,.$ If $p\in(-1,1)$, then
$$a_1^p+a_2^p+\ldots+a_n^p\leq 
b_1^p+b_2^p+\ldots+b_n^p\,.$$
This inequality fails for $-1<p<0$ (but remains true
for $0<p<1$) if the constraint 
$a_1a_2\ldots a_n=b_1b_2\ldots b_n$ is replaced by
the weaker one $a_1a_2\ldots a_n\leq b_1b_2\ldots b_n\,.$
\end{cor}

\noindent {\bf Proof of  \theref{theorem_1} } 
If $S$ is empty, then $E_k(a)=E_k(b)$ for all $k\in\{0,1,\dotsc,n-1\}$ and hence $a=b$, which immediately implies the inequality. We therefore assume that $S$ is nonempty.

Let $y\colon\,[0,1]\to \Delta_n$ be the curve connecting points $a$ and $b$
as in Definition \ref{definition_1}. Consider the function
\begin{align}
   p(t,x)&=(x+y_1(t))(x+y_2(t))\ldots (x+y_n(t))=\sum_{k=0}^{n-1} x^kE_k(y(t)) + x^n  \notag \\
             &=(x+a_1)(x+a_2)\ldots(x+a_n)+\sum_{k\in S} x^k A_k(t)\,,
\end{align}
where $A_k(t)=E_k(y(t))-E_k(a)$ is a non-decreasing mapping. Our goal is to show that the
function 
\begin{align}
\eta(t)=\sum_{i=1}^n f(y_i(t))
\end{align}
is non-decreasing on $[0,1]$, i.e. we show that $\eta'(t)\geq 0$ a.e. on $(0,1)$.

To this end, fix $i\in\{1,2,\ldots,n\}$. Since $p(t,-y_i(t))=0$ for all $t\in(0,1)$, 
we obtain
$$\partial_1\,p(t,-y_i(t))+\partial_2\,p(t,-y_i(t))\cdot (-y_i'(t))=0$$
for all $t\in(0,1)$ and therefore
\begin{align}
 \sum_{k\in S} (-y_i(t))^k A'_k(t)+\prod_{j\neq i} (y_j(t)-y_i(t)) \cdot (-y_i'(t))=0\,,
 \end{align}
which gives
$$y_i'(t)=\sum_{k\in S} (-y_i(t))^k A'_k(t)\Bigl(\prod_{j\neq i} (y_j(t)-y_i(t))\Bigr)^{-1}\,.$$
This equality holds, if $y_i(t)\neq y_j(t)$ for $i\neq j$, which
is true for all but countably many values of $t\in(0,1)$.
For those values of $t$ we get
\begin{align}
    \eta'(t)&=\sum_{i=1}^n f'(y_i(t))\cdot y_i'(t) \notag\\
        &=\sum_{i=1}^n f'(y_i(t))\cdot\sum_{k\in S} (-y_i(t))^k A'_k(t)\Bigl(\prod_{j\neq i} (y_j(t)-y_i(t))\Bigr)^{-1}\notag\\
      &=\sum_{k\in S} A'_k(t)\sum_{i=1}^n f'(y_i(t))\cdot (-y_i(t))^k\, \Bigl(\prod_{j\neq i} (y_j(t)-y_i(t))\Bigr)^{-1}\,.
   \end{align}   
Fix $t\in(0,1)$ such that $y_i(t)\neq y_j(t)$ for $i\neq j$
and write $y_i=y_i(t)$ for simplicity. Since $A'_k(t)\geq 0$, we will be done, if we show that
$$\widehat{D}:=\sum_{i=1}^n f'(y_i)\cdot (-y_i)^k \Bigl(\prod_{j\neq i} (y_j-y_i)\Bigr)^{-1}\geq 0\quad\hbox{for all $k\in S$}\,.$$

To this end, consider the polynomial
$$g(x)=\sum_{i=1}^n f'(y_i)\cdot (-y_i)^k \Bigl(\prod_{j\neq i} (y_j-y_i)\Bigr)^{-1}\cdot \prod_{j\neq i} (x-y_j)\,.$$
The degree of $g$ equals $n{-}1$ and the coefficient at $x^{n-1}$ is equal to $\widehat{D}$. Moreover,
$$g(y_i)= f'(y_i)\cdot (-y_i)^k\cdot (-1)^{n-1}\quad (i=1,2,\ldots,n)\,.$$
Therefore the function $h(x)=g(x)+(-1)^{n+k}x^kf'(x)$ has $n$ different roots $y_1,y_2,\ldots,y_n$
in the interval~$I$.
It follows that the function 
\begin{align}
    h^{(n-1)}(x) =(n-1)! \, \widehat{D} + (-1)^{n+k}(x^kf'(x))^{(n-1)}
\end{align}    
has a root in the
interval $I$, and since $(-1)^{n+k}(x^kf'(x))^{(n-1)}\leq 0$ for all $x\in I$, 
it follows that $\widehat{D}\geq 0$, which completes the proof of \theref{theorem_1}. \qed

\noindent {\bf Proof of \theref{theorem_2} } 
Suppose, to the contrary, that $(-1)^{k+n}(x^kf'(x))^{(n-1)}>0$ for some $x\in I$ and some $k\in S$.
Then $(-1)^{k+n}(x^kf'(x))^{(n-1)}>0$ holds for all $x$ belonging to some interval $J$ contained in~$I$.
Choose the numbers $a_1<a_2<\ldots<a_n$ from $J$ and consider
$$p(t,x)=(x+a_1)\cdot(x+a_2)\cdot\ldots\cdot(x+a_n)+t\,x^k\,.$$
Then for all sufficiently small $t$ $(0<t<\varepsilon)$, there exist different numbers
$y_i(t)$ belonging to~$J$, such that
$$p(t,x)=(x+y_1(t))(x+y_2(t))\ldots(x+y_n(t))\,.$$
Then
\[
	x^n + \sum_{i=0}^{n-1} E_i(a) \cdot x^i + t\,x^k = p(t,x) = x^n + \sum_{i=0}^{n-1} E_i(y(t)) \cdot x^i\,,
\]
and since $t>0$, we see that $a$ and $b=y(t)$ satisfy \eqref{dupa}.
We will be done if we show that
$$f(a_1)+f(a_2)+\ldots+f(a_n)>f(y_1(t))+f(y_2(t))+\ldots+f(y_n(t))\,.$$

We proceed in the same way as in the proof of \theref{theorem_1}. We define
$$\eta(t)=\sum_{i=1}^n f(y_i(t))\quad\hbox{for $0<t<\varepsilon$}\,$$
and this time we want to show that $\eta'(t)<0$ for $0<t<\varepsilon$.

By the Inverse Mapping Theorem (see proof of \propref{proposition_2} below for a more
detailed explanation), $y\in C^1(0,\varepsilon)$ and therefore
\begin{align}
   \eta'(t)=\sum_{i=1}^n f'(y_i(t))\cdot y_i'(t)=\sum_{i=1}^n f'(y_i(t))\cdot (-y_i(t))^k 
\Bigl(\prod_{j\neq i} (y_j(t)-y_i(t))\Bigr)^{-1}\,.
\end{align}
Now, like previously, write $y_i=y_i(t)$ for simplicity. Our goal is therefore to prove that
$$\widehat{D}:=\sum_{i=1}^n f'(y_i)\cdot (-y_i)^k \Bigl(\prod_{j\neq i} (y_j-y_i)\Bigr)^{-1}< 0\,.$$
Consider the polynomial
$$g(x)=\sum_{i=1}^n f'(y_i)\cdot (-y_i)^k \Bigl(\prod_{j\neq i} (y_j-y_i)\Bigr)^{-1}\cdot \prod_{j\neq i} (x-y_j)\,.$$
The degree of $g$ equals $n{-}1$ and the coefficient at $x^{n-1}$ is equal to $\widehat{D}$. Moreover,
the function $h(x)=g(x)+(-1)^{n+k}x^kf'(x)$ has $n$ different roots $y_1,y_2,\ldots,y_n$
in the interval~$J$.
It follows that the function 
$$h^{(n-1)}(x) =(n-1)!\, \widehat{D} + (-1)^{n+k}(x^kf'(x))^{(n-1)}$$ has a root in the
interval $J$. And since $(-1)^{n+k}(x^kf'(x))^{(n-1)}>0$ for all $x\in J$, 
it follows that $\widehat{D}<0$, which completes the proof of \theref{theorem_2}.    \qed

\section{Construction of the connecting curve}

In this section we prove that
condition \eqref{inequality_1} implies $a\preceq b$, if $2\leq n\leq 4$,
$I=(0,\infty)$ and $S\subseteq\{1,2,\ldots,n{-}1\}$.
However, we start with a construction of the desired
curve for a general interval $I$, integer $n\geq 2$
and set $S\subseteq\{0,1,\ldots,n{-}1\}$.

For $a,b\in \Delta_n$, we say that $a<b$, if $a\neq b$ and 
$E_k(a)\leq E_k(b)$ for all $k=0,1,\ldots,n{-}1$. We say that $a\leq b$,
if $a<b$ or $a=b$.

\begin{defi}
\label{definition_2}
For $a<b$ denote by $\C(a,b)$ the set of all piecewise differentiable 
(i.e. continuous and differentiable in all but at most countably many points) curves
$y$ in $\Delta_n$ satisfying:

(a) the curve $y(t)$ starts at $a$ (i.e. $y(0)=a$, if the curve $y(t)$ is parametrized by the
interval $[0,\ve]$);

(b) $y(t)\in\inter(\Delta_n)$ for all but at most countable many values $t$;

(c) the mappings $E_k(y(t))$ are non-decreasing in $t$
and $E_k(y(t))\leq E_k(b)$ for all $t$ and each $k=0,1,\ldots,n{-}1$.

\noindent Note that a curve in $\C(a,b)$ does not necessarily end at the point $b$.
\end{defi}

\begin{prop}
\label{proposition_1}
Let $n\geq 2$ be a positive integer and let
$S$ be a nonempty subset of $\{0,1,\ldots,n{-}1\}$.
Let moreover $a,b\in \Delta_n$ be such that \eqref{inequality_1} holds.
Furthermore, suppose that for all $c\in \Delta_n$
with $a\leq c<b$ the set $\C(c,b)$ is nonempty. 
Then $a\preceq b$.
\end{prop}

\noindent {\bf Proof.} 
Each element (curve) of $\C(a,b)$ is a (closed) subset of $\Delta_n$.
We equip the set $\C(a,b)$ with the inclusion relation $\subseteq$,
obtaining a nonempty partially ordered set $(\C(a,b),\subseteq)$. 
We are going to show that 
each chain $\{y_i\}_{i\in \mathcal{I}}$ has an upper bound in $\C(a,b)$.

To achieve this, consider the curve
$$y_0=\overline{\bigcup_{i\in \mathcal{I}} y_i}\,,$$
i.e. the concatenation of the curves $y_i$. Then obviously $y_0$ satisfies
conditions (a) and (c) of Definition \ref{definition_2}.
To prove (b) assume that $y_0$ is parametrized on $[0,1]$. Then for each positive
integer $k$ the curve $y_k$, defined as the  restriction of $y_0$ 
to the interval $[0,1-{1\over k}]$,
is contained in some curve $y_i\in \C(a,b)$ of the given chain $\{y_i\}$. Therefore
$y_k(t)$ is piecewise differentiable and satisfies condition (b) for each
positive integer $k$. Moreover,
$$y_0=\overline{\bigcup_{k=1}^\infty y_k}\,.$$
Hence $y_0$ is piecewise differentiable and satisfies (b) as well.

Now, by the Kuratowski-Zorn lemma,
there exists a maximal element $y$ in $(\C(a,b),\subseteq)$. We show that $y$
is a desired curve connecting the points $a$ and $b$, which will imply
that $a\preceq b$.

To this end, it is enough to show that, if the curve $y$ is parametrized on $[0,1]$,
then $y(1)=b$. Suppose, to the contrary, that $y(1)=c\neq b$. Then $a\leq c<b$,
and hence the set $\C(c,b)$ is nonempty. Thus the curve $y$ can be extended beyond
the point $c$, which contradicts the fact that $y$ is a maximal element in $\C(a,b)$.
This completes the proof of \propref{proposition_1}.    \qed

\noindent From now on assume that $I=(0,\infty)$ and $S$ is a nonempty subset of
$\{1,2,\ldots,n{-}1\}$.

In order to prove that \eqref{inequality_1} implies $a\preceq b$,
it suffices to show that the sets $\C(a,b)$ for $a,b\in\Delta_n$
with $a<b$ are nonempty. This is implied by the following conjecture,
which we will prove later for $n\leq 4$.

\begin{conj}
\label{main_conjecture}
Let $n\geq 2$ be an integer and $a\in\Delta_n$.
Let $S$ be a nonempty subset of $\{1,2,\ldots,n{-}1\}$
with the property that there exist $A_k>0$ for $k\in S$
such that all the roots of the polynomial 
$$q(x)=(x+a_1)(x+a_2)\ldots(x+a_n)+\sum_{k\in S} A_kx^k$$
are real (and hence negative). Then there exist 
continuous on $[0,\ve]$, differentiable on $(0,\ve)$
and nondecreasing mappings $B_k:[0,\ve]\to\R$ ($k\in S$)
with $B_k(0)=0$ such that
$\sum_{k\in S} B_k(t)$ is increasing on $[0,\ve]$
and for all sufficiently small values of $t>0$ the polynomial
$$(x+a_1)(x+a_2)\ldots(x+a_n)+\sum_{k\in S} B_k(t)x^k$$
has $n$ distinct real (and hence negative) roots.
\end{conj}

\smallskip
Now we show how \conjref{main_conjecture} implies
that the sets $\C(a,b)$ are nonempty.

\begin{prop}
\label{proposition_2}
Let $n$ and $S$ be such that the conjecture holds. 
Let moreover $a,b\in \Delta_n$ be such that \eqref{inequality_1} holds.
Then the set $\C(a,b)$ is nonempty.
\end{prop}
\noindent {\bf Proof.} 
Consider the polynomials 
$$p(x)=(x+a_1)(x+a_2)\ldots(x+a_n)\quad\hbox{and}\quad
q(x)=(x+b_1)(x+b_2)\ldots(x+b_n)\,.$$
Then
$$q(x)-p(x)=\sum_{k=0}^{n-1} (E_k(b)-E_k(a))x^k=\sum_{k\in S} A_kx^k\,,$$
where $A_k>0$ for all $k\in S$.
According to the conjecture, there exist
continuous on $[0,\ve]$ and differentiable on $(0,\ve)$
nondecreasing mappings $B_k:[0,\ve]\to\R$, 
with $B_k(0)=0$ such that $\sum_{k\in S}B_k(t)$ is increasing on $[0,\ve]$
and for all $t\in(0,\ve)$ the polynomial
$$p(x)+ \sum_{k\in S} B_k(t)x^k$$
has $n$ distinct real (and hence negative) roots $-y_n(t)<-y_{n-1}(t)<\ldots<-y_1(t)<0$.
We show that $y(t)=(y_1(t),y_2(t),\ldots,y_n(t))$ defines a
differentiable curve (parametrized on $[0,\ve]$) that belongs to $\C(a,b)$,
provided $\ve$ is chosen in such a way that $B_k(\ve)\leq A_k$ for $k\in S$.

Consider the mapping $\Psi\colon\,\overline{\Delta_n}\to \Psi(\overline{\Delta_n})$ 
given by
$$\Psi(y)=(E_{n-1}(y),E_{n-2}(y),\ldots,E_0(y))\,.$$
Then it follows from Remark \ref{injectivity_remark} that the mapping $\Psi$ is injective, hence $\Psi$ is a continuous bijection defined on a closed subset of $\R^n$.
Therefore the restriction $\left.\Psi\right|_U$ of $\Psi$ to a neighbourhood $U$ of $a$ is continuously invertible and thus
$$y(t)=\Psi^{-1}(\Psi(a)+(B_0(t),B_1(t),\ldots,B_{n-1}(t)))\quad\hbox{($t\in[0,\ve]$)}$$
(here we put $B_k(t)=0$ for $k\not\in S$)
is a curve starting at $a$; note that $\Psi(a)+(B_0(t),B_1(t),\ldots,B_{n-1}(t))$ is contained in $\Psi(U)$ for sufficiently small $\ve$.
Moreover $y(t)\in\Delta_n$. 
Hence condition (a) is satisfied.
Since $y(t)\in \inter(\Delta_n)$ for all $t\in(0,\ve)$,
condition (b) holds. It is also clear that (c) is satisfied, since
$E_k(y(t))=E_k(a)+B_k(t)\leq E_k(a)+A_k=E_k(b)$ for all $k\in\{0,1,\ldots,n{-}1\}$.

It remains to prove that $y(t)$ is differentiable on $(0,\ve)$. This however
is a consequence of the Inverse Mapping Theorem, if we show that
$$\det{D\Psi(y)}\neq 0 \quad\hbox{for all $y\in \inter(\Delta_n)$.}$$

To this end, 
let $V(y)$ be the $n\times n$ Vandermonde-type matrix given by 
$V_{ij}(y)=(-y_i)^{n-j}$ ($1\leq i,j\leq n$). This matrix is obtained
from the standard Vandermonde matrix
\begin{align}
W (-y_1, -y_2, \ldots , -y_n) =
\begin{pmatrix}
  1      & -y_1    & (-y_1)^2  & \cdots & (-y_1)^{n-1} \\
  1      & -y_2    & (-y_2)^2  & \cdots & (-y_2)^{n-1} \\
  1      & -y_3    & (-y_3)^2  & \cdots & (-y_3)^{n-1} \\
  \vdots & \vdots & \vdots & \ddots & \vdots    \\
  1      & -y_n    & (-y_n)^2  & \cdots & (-y_n)^{n-1}
\end{pmatrix}
\end{align}
by reversing the order of columns of $W$.

\newcommand{\pdd}[1]{\frac{\partial}{\partial #1}}%
Since \cite{dannan2015}
\[
	(D\Psi(y))_{jk} = \pdd{y_k}\,E_{n-j}(y) =
	\begin{cases}
		1 &: j=1\,,\\
		E_{n-j}(y^{(k)}) &: j>1\,,
	\end{cases}
\]
where $y^{(k)} = (y_1,\dotsc,y_{k-1},y_{k+1},\dotsc,y_n)$ is $y$ with its $k$-th component removed, it follows from the general formula

\begin{align}
    t^{n-1}+\sum_{j=0}^{n-2}t^jE_j(z_1,z_2,\ldots,z_{n-1})=(t+z_1)(t+z_2)\ldots(t+z_{n-1})
 \end{align}   
that
\begin{align*}
	(V(y) \cdot D\Psi(y))_{ik} &= \sum_{j=1}^n (V(y))_{ij} \cdot (D\Psi(y))_{jk}\\
	&= (-y_i)^{n-1} + \sum_{j=2}^n (-y_i)^{n-j} \cdot E_{n-j}(y^{(k)})\\
	&= (-y_i)^{n-1} + \sum_{j=0}^{n-2} (-y_i)^{j} \cdot E_{j}(y^{(k)})%
	= \prod_{j\neq k} (y_j-y_i)
\end{align*}
and thus

\begin{align}
    V(y)\cdot D\Psi(y)=
\diag\Bigl(\prod_{j\neq 1}(y_j-y_1),\prod_{j\neq 2}(y_j-y_2),\ldots,\prod_{j\neq n}(y_j-y_n)\Bigr)\,.
\end{align}
It is well-known that
$$\det{V(y)}=\prod_{i<j}(y_j-y_i)\neq 0 \quad\hbox{($y\in\inter{\Delta_n}$)}\,.$$
Therefore we obtain
$$\det {D\Psi(y)}=\prod_{i<j}(y_i-y_j)\neq 0\quad\hbox{($y\in\inter{\Delta_n}$)}\,,$$
which completes the proof of \propref{proposition_2}. \qed

\begin{lem}
\label{lemma_1}
Assume that $n\geq 3$ is odd and let $0<a_1\leq a_2\leq\ldots\leq a_n$.
Let moreover $A_k\geq 0$ for $k=1,2,\ldots,(n{-}1)/2$ with at least one $A_k$ not equal to $0$. Consider the polynomials
\begin{align}
   P(x)&=(x+a_1)(x+a_2)\ldots(x+a_n)+\sum_{k=1}^{(n-1)/2} A_kx^{2k-1}\,,\notag\\
   Q(x)&=(x+a_1)(x+a_2)\ldots(x+a_n)+\sum_{k=1}^{(n-1)/2} A_kx^{2k}\,.
\end{align}
Then the polynomial $P$ has exactly one root in the interval $(-a_1,0)$ and at most
two roots in  the interval $(-a_n,-a_{n-1})$. Moreover, the polynomial $Q$ has
exactly one root in the interval $(-\infty,-a_n)$ and at most two 
roots in the interval $(-a_2,-a_1)$.
\end{lem}

\noindent {\bf Proof.} 
That $P$ has exactly one root in $(-a_1,0)$ follows immediately from the observation
that $P(-a_1)<0$, $P(0)>0$ and $P'(x)>0$ on $(-a_1,0)$. 

Now we show that $Q$ has exactly one root in $(-\infty,-a_n)$. 

Dividing the
equation $Q(x)=0$ by $x^na_1a_2\ldots a_n$ and substituting $z=1/x$ and $b_i=1/a_i$, yields
the equation $P_0(z)=0$, where
$$P_0(z)=(z+b_1)(z+b_2)\ldots(z+b_n)+\sum_{k=1}^{(n-1)/2} B_kz^{2k-1}$$
for some nonnegative numbers $B_k$, not all equal to $0$.
We already know that $P_0$ has exactly one root in the interval $(-b_n,0)$,
so it follows that $Q$ has exactly one root in the interval $(-\infty,-a_n)$.

Now we prove that $Q$ has at most two roots in the interval $(-a_2,-a_1)$.
To the contrary, suppose that $Q$ has at least $3$ roots in $(-a_2,-a_1)$.
Since $Q(-a_2)>0$ and $Q(-a_1)>0$, it follows that $Q$ has an even
number, and hence at least four, roots in the interval $(-a_2,-a_1)$.

Let $0>-c_1\geq-c_2\geq\ldots\geq-c_{n-1}$ be the roots of $p'(x)=0$, where
\begin{align}
    p(x)=(x+a_1)(x+a_2)\ldots(x+a_n)\,.
 \end{align}   
Then $a_1<c_1<a_2$. The polynomial $Q(x)$ is decreasing on the interval 
$[-a_2,-c_1]$, so it has at most one root in this interval. Therefore
the polynomial $Q$ has at least three roots in the interval $(-c_1,-a_1)$,
and consequently the equation $Q''(x)=0$ has a root in $(-c_1,-a_1)$.
But $Q''(x)>0$ for all $x>-c_1$, a contradiction. Hence $Q$ must have at most
two roots in $(-a_2,-a_1)$.

Finally, to prove that $P$ has at most two roots in the interval 
$(-a_n,-a_{n-1})$, divide the
equation $P(x)=0$ by $x^na_1a_2\ldots a_n$ and substitute $z=1/x$ and $b_i=1/a_i$.
This reduces to the equation $Q_0(z)=0$, where
$$Q_0(z)=(z+b_1)(z+b_2)\ldots(z+b_n)+\sum_{k=1}^{(n-1)/2} B_kz^{2k}$$
for some nonnegative numbers $B_k$, not all equal to $0$.
We already know that $Q_0$ has at most two roots in the interval $(-b_{n-1},-b_n)$,
so it follows that $P$ has at most two roots in the interval $(-a_n,-a_{n-1})$.
This completes the proof of \leref{lemma_1}.  \qed

\smallskip
The same proof yields an analogous result for even values of $n$.

\begin{lem}
\label{lemma_2}
Assume that $n\geq 2$ is even and let $0<a_1\leq a_2\leq\ldots\leq a_n$.
Let moreover $A_k\geq 0$ for $k=1,2,\ldots,n/2$ and not all of the $A_k$'s
are equal to $0$. Consider the polynomials
\begin{align}
    P(x)&=(x+a_1)(x+a_2)\ldots(x+a_n)+\sum_{k=1}^{n/2} A_kx^{2k-1}\,,\notag\\
     Q(x)&=(x+a_1)(x+a_2)\ldots(x+a_n)+\sum_{k=1}^{n/2-1} A_kx^{2k}\,.
\end{align}
Then the polynomial $P$ has exactly one root in each of the intervals 
$(-\infty,-a_n)$ and $(-a_1,0)$ and $Q$ has at most
two roots in each of the intervals $(-a_n,-a_{n-1})$ and $(-a_2,-a_1)$.
\end{lem}
{\bf Proof.} The same proof as that for \leref{lemma_1} can be used. \qed 

Now we turn to the proof of \conjref{main_conjecture} for $2\leq n\leq 4$ and
an arbitrary nonempty set $S\subseteq\{1,2,\ldots,n{-}1\}$.

We first make some useful general remarks.

Let $I(a)=\{i\in\{1,2,\ldots,n{-}1\}\;:\;a_i=a_{i+1}\}$.
If $I(a)$ is empty, then the conjecture holds. Indeed, if $k\in S$,
then all the roots of the polynomial
$$(x+a_1)(x+a_2)\ldots(x+a_k)+t\,x^k$$
are, for all sufficiently small $t>0$, real and distinct.

\smallskip
On the other hand, if $I(a)=\{1,2,\ldots,n{-}1\}$, then only the set
$S=\{1,2,\ldots,n{-}1\}$ possibly satisfies the assumptions of the conjecture. 
Indeed, suppose that $l\not\in S$ and let
$-b_1\geq-b_2\geq\ldots\geq-b_n$ be the roots of 
$$q(x)=(x+a_1)^n+\sum_{k\in S} A_kx^k\,.$$
Then by the inequality of arithmetic and geometric means, we obtain
\begin{align}
   {E_l(a)\over{n\choose l}}={E_l(b)\over{n\choose l}}
\geq (E_0(b))^{(n-l)/n}=(E_0(a))^{(n-l)/n}={E_l(a)\over{n\choose l}}\,,
\end{align}
and hence $b_1=b_2=\ldots=b_n$. Since $E_0(a)=E_0(b)$, it follows that $a=b$, i.e.
$A_k=0$ for all $k\in S$. A contradiction.

\smallskip
Let $I$ be a non-empty subset of $\{1,2,\ldots,n{-}1\}$.
We observe that the conjecture is true for a set $S$ and all 
$a\in\Delta_n$ with $I(a)=I$, if it is true for
a set $T=\{n{-}k\;:\;k\in S\}$ and all $b\in\Delta_n$
with $I(b)=\{n{-}i\;:\;i\in I\}$. Indeed: if all the roots of the polynomial
$$q(x)=(x+a_1)(x+a_2)\ldots(x+a_n)+\sum_{k\in S} A_kx^k$$
are real, then substituting $x=1/z$ and $a_i=1/b_i$, we infer that
all the roots of the polynomial
$$r(z)=(z+b_1)(z+b_2)\ldots(z+b_n)+\sum_{l\in T} B_lz^l$$
are real. Hence there exist continuous on $[0,\ve]$, differentiable
on $(0,\ve)$ and nondecreasing 
mappings $C_l(t)$ with $C_l(0)=0$ such that
the polynomial 
$$(z+b_1)(z+b_2)\ldots(z+b_n)+\sum_{l\in T} C_l(t)z^l$$
has $n$ distinct real roots. Substituting $z=1/x$ and $b_i=1/a_i$,
we infer that the polynomial
$$(x+a_1)(x+a_2)\ldots(x+a_n)+\sum_{k\in S} C_{n-k}(t)x^k$$
has $n$ distinct real roots.

\smallskip
For $n=2$ the only possibility for the set $S$ is $\{1\}$ and it is enough to
notice that the polynomial $(x+a_1)(x+a_2)+t\,x$ has two distinct real roots for any
$t>0$.

\medskip
Assume now $n=3$. Then, in view of the above remarks, we have to consider 
two cases:
1) $a_1<a_2=a_3$;
2) $a_1=a_2=a_3$.

\smallskip

1) If $2\notin S$, then the condition of \conjref{main_conjecture} can not 
be satisfied since for $A_1>0$, according \leref{lemma_1}, the polynomial
$$P(x) = (x+a_1)(x+a_2)^2+A_1x$$
has only one real root in the interval $(-a_1,0)$ and
obviously no roots on $\R\setminus(-a_1,0)$. Thus $P$ has
only one real root for all $A_1>0$. We can therefore assume $2\in S$, and for all sufficiently small $t>0$, the polynomial
$$(x+a_1)(x+a_2)^2+t\,x^2$$
has three distinct real roots.

\smallskip
2) According to the above remarks, $S=\{1,2\}$. Then the polynomial
$(x+a_1)^3+t\,a_1x+t\,x^2$ has 3 distinct real roots for all sufficiently small
$t>0$.

\medskip
Assume $n=4$. In this case we have 5 possibilities:
1)~$a_1=a_2<a_3<a_4$;
2)~$a_1<a_2=a_3<a_4$;
3)~$a_1<a_2=a_3=a_4$;
4)~$a_1=a_2<a_3=a_4$;
5)~$a_1=a_2=a_3=a_4$.

\smallskip
1) We note that $S\neq\{2\}$, since, by \leref{lemma_2}, the polynomial
$$Q(x)=(x+a_1)^2(x+a_3)(x+a_4)+A_2x^2\quad\hbox{for $A_2>0$}$$
has at most two real roots in the interval $(-a_4,-a_3)$ and obviously
no roots on $\R\setminus(-a_4,-a_3)$. Thus $Q$ has at most two real roots. 
Therefore $S$ contains an odd integer $k$. Then for all 
sufficiently small $t>0$, the polynomial $(x+a_1)^2(x+a_3)(x+a_4)+t\,x^k$
has four distinct real roots. 

\smallskip
2) Note that $2\in S$, since by \leref{lemma_2}, the polynomial
$$(x+a_1)(x+a_2)^2(x+a_4)+A_1x+A_3x^3\quad\hbox{for $A_1,A_3>0$}$$
has at most two real roots. Then for all 
sufficiently small $t>0$, the polynomial
$$(x+a_1)(x+a_2)^2(x+a_4)+t\,x^2$$
has four distinct real roots. 

\smallskip
3) We observe that $\{1,2\}\subset S$ or $\{2,3\}\subset S$, 
since by \leref{lemma_2}, each of the polynomials
$$(x+a_1)(x+a_2)^3+A_1x+A_3x^3\quad\hbox{and}\quad
(x+a_1)(x+a_2)^3+A_2x^2 \quad\hbox{for $A_1,A_2,A_3>0$}$$
as well as
$$(x+a_1)(x+a_2)^3+A_1x\quad\hbox{and}\quad
(x+a_1)(x+a_2)^3+A_3x^3 \quad\hbox{for $A_1,A_3>0$}$$
has at most two real roots. Moreover, we prove that $S\neq \{1,2\}$.

Suppose that the polynomial
$Q(x)=(x+a_1)(x+a_2)^3+A_1x+A_2x^2$
has four real roots. Let $Q_1(x)=(x+a_1)(x+a_2)^3$ and $Q_2(x)=A_1x+A_2x^2$.
Let $-c\neq a_2$ be the root of the polynomial $Q_1'(x)$ and
let $-d$ be the root of $Q_2'(x)$.

If $d<c$, then $Q$ is decreasing on $(-\infty,-c]$, so $Q$ has at most
one root in this interval. Therefore $Q$ has at least 3 roots in the interval
$(-c,0)$. Thus $Q''(x)$ has a root in the interval $(-c,0)$, which is impossible,
since $Q''(x)>0$ on $(-c,0)$.

If $a_2\geq d\geq c$, then $Q$ is increasing on the interval
$[-c,0)$ and decreasing on the interval $(-\infty,-d]$, so $Q$ must have at least two roots
in the interval $(-d,-c)$. But $Q(x)<0$ on this interval.

Finally, if $d>a_2$, then $Q$ may only have roots in the union
$(-\infty,a_2)\cup(-a_1,0)$. But $Q$ is increasing on $(-a_1,0)$,
so $Q$ has 3 roots in $(-\infty,a_2)$. This however is impossible, since
$Q''(x)>0$ for $x\in(-\infty,a_2)$. Thus $\{2,3\}\subseteq S$
and the polynomial
$$(x+a_1)(x+a_2)^3+t\,x^2(x+a_2)$$
has, for all sufficiently small $t>0$, four distinct roots.

\smallskip
4) Since the polynomial
$(x+a_1)^2(x+a_3)^2+A_2x^2$ has no real roots,
$1\in S$ or $3\in S$. Then the polynomial
$(x+a_1)^2(x+a_3)^2+t\,x^k$ for $k=1,3$
has, for all sufficiently small $t>0$, four distinct real roots.

\smallskip
5) In view of the above remarks, $S=\{1,2,3\}$. Consider
$$r(x)=(x+a_1)^4+t\,x^3+2\,t\,a_1x^2+t(a_1^2-t^2)x=(x+a_1)^4+t\,x((x+a_1)^2-t^2)\,.$$
Then for all sufficiently small $t>0$, $a_1^2-t^2>0$, and the polynomial
$r$ has four distinct real roots, because
$$r(-a_1-2t)=t^3(10t-3a_1)<0\,,\ 
r(-a_1)=a_1t^3>0\ \hbox{and}\ 
r(-a_1+2t)=t^3(22\,t-3a_1)<0\,.$$

Thus we have proved:

\begin{cor}
\label{corollary_3}
\conjref{main_conjecture}
 is true if $2\leq n\leq 4$ and $S$ is an arbitrary nonempty
subset of $\{1,2,\ldots,n{-}1\}$.
\end{cor}
This implies that the sum of squared logarithms inequality (\conjref{log_inequality}) holds also for $n=4$.
\begin{cor}[Sum of squared logarithms inequality for $n=4$]
Let $a_1,a_2,a_3,a_4,b_1,b_2,b_3,b_4>0$ be given positive numbers such that
\begin{align}
   a_1+a_2+a_3+a_4 &\le b_1+b_2+b_3+b_4  \, , \notag\\
   a_1\, a_2+a_1\, a_3+ a_2\, a_3+a_1\, a_4+a_2\, a_4+a_3\, a_4  & \le b_1\, b_2+b_1\, b_3+ b_2\, b_3+b_1\, b_4+b_2\, b_4+b_3\, b_4 \, ,\notag \\
    a_1 \, a_2\, a_3+a_1\, a_2\, a_4+a_2\, a_3\, a_4+a_1\, a_3\, a_4 &\le b_1\, b_2\, b_3+b_1\, b_2\, b_4+b_2\, b_3\, b_4+b_1\, b_3\, b_4  \, , \notag\\
a_1 \, a_2\, a_3\, a_4&= b_1\, b_2\, b_3\, b_4\, . \notag
\end{align}
Then
\begin{align}
 \log^2 a_1+\log^2 a_2+\log^2 a_3+\log^2 a_4     \le  \log^2 b_1+\log^2 b_2+\log^2 b_3+\log^2 b_4   \, .\notag
 \end{align}
\end{cor}
\noindent {\bf Proof.} Use \coref{corollary_3} and observe that $S$ may be an arbitrary subset of $\{1,2,3\}$. \qed

\smallskip

\begin{cor}
\label{corollary_4}
Let $n\geq 2$ be an integer and let $T$ be an arbitrary subset of
$\{1,2,\ldots,n{-}1\}$. Assume that the Conjecture \ref{main_conjecture} holds
for $n$ and for any nonempty subset $S$ of $T$. Let moreover $f\in C^n(0,\infty)$.
Then the inequality
$$f(a_1)+f(a_2)+\ldots+f(a_n)\leq f(b_1)+f(b_2)+\ldots+f(b_n)$$
holds for all $a,b\in\Delta_n$ satisfying
\begin{align}
\label{inequality_4}
   E_k(a)\leq E_k(b) \quad\hbox{for } k\in T\quad \hbox{and} \quad E_k(a)=E_k(b) 
\quad\hbox{for } k=0 \hbox{ or }k\not\in T
\end{align}
if and only if
\begin{align}
 \label{inequality_5}
   (-1)^{n+k}(x^k f'(x))^{(n-1)}\leq 0 \quad \hbox{for all }x>0 \hbox{ and all } k\in T\,.
\end{align}
\end{cor}

\noindent {\bf Proof.} 
Assume first \eqref{inequality_5} holds and let $a,b\in\Delta_n$
satisfy \eqref{inequality_4}. Consider any $c\in\Delta_n$ with $a\leq c<b$.
Then the pair $c$, $b$ satisfies condition \eqref{inequality_1} for some nonempty
subset $S$ of $T$. Therefore by \propref{proposition_2}, the set
$\C(c,b)$ is nonempty and hence by \propref{proposition_1}, $a\preceq b$.
Now \theref{theorem_1} implies that inequality \eqref{inequality_3} holds.

Conversely, if \eqref{inequality_3} holds for all $a,b\in\Delta_n$
satisfying \eqref{inequality_4}, then \eqref{inequality_3} also holds for all $a,b\in\Delta_n$
satisfying condition \eqref{inequality_1} with $S=T$. Thus \theref{theorem_2} implies \eqref{inequality_5}. This completes
the proof.   \qed

\section{Outlook}
Our result generalizes and extents the previously 
known results on the sum of squared logarithms inequality. 
Indeed, compared to the proof in \cite{Neff_log_inequality13} 
our development here views the problem from a different angle 
in that it is not the logarithm function that defines the problem, 
but a certain monotonicity property in the geometry of
polynomials, explicitly stated in \conjref{main_conjecture}. 

If one tries to adopt the above proof of
\conjref{main_conjecture} for $n\leq 4$ to the case $n\geq 5$, one has to
deal with approximately $2^n$ cases considered separately.
Therefore it is clear, that the extension to natural numbers 
$n$ beyond $n=6$, say, is out of reach with such a method. 
Instead, a general argument should be found to
prove or disprove \conjref{main_conjecture} for general $n$.
Furthermore, it might be worthwhile to develop a better understanding of the differential inequality condition $(-1)^{n+k}(x^k f'(x))^{(n-1)}\leq 0$.

\section*{Competing interests}
The authors declare that they have no competing interests.

\section*{Authors' contributions}
Both authors contributed fully to all parts of this paper.

\section*{Acknowledgements}
We thank Johannes Lankeit (Universität Paderborn) as well as Robert Martin (Universität Duisburg-Essen) for their help in revising this paper.

\bibliographystyle{plain} %
{ %

}

\end{document}